\documentclass{IEEEtran}

\usepackage{graphicx}
\usepackage{amssymb}
\usepackage{amsmath}
\usepackage{subfigure}
\usepackage{epsfig}
\usepackage{color}
\usepackage{enumitem}
\usepackage{float}
\usepackage{soul}
\usepackage{wrapfig}
\usepackage{arydshln}
\usepackage{url}

\newtheorem{Remark 1}{Remark}
\newtheorem{Remark 2}[Remark 1]{Remark}
\newtheorem{Remark 3}[Remark 1]{Remark}
\newtheorem{Remark 4}[Remark 1]{Remark}
\newtheorem{Remark 5}[Remark 1]{Remark}
\newtheorem{Remark 6}[Remark 1]{Remark}
\newtheorem{Remark 7}[Remark 1]{Remark}

\newtheorem{Assumption 1}{Assumption}
\newtheorem{Definition 1}{Definition}
\newtheorem{Theorem 1}{Theorem}
\newtheorem{Theorem 2}[Theorem 1]{Theorem}
\newtheorem{Theorem 3}[Theorem 1]{Theorem}
\newtheorem{Theorem 4}[Theorem 1]{Theorem}
\newtheorem{Theorem 5}[Theorem 1]{Theorem}
\newtheorem{Theorem 6}[Theorem 1]{Theorem}
\newtheorem{Theorem 7}[Theorem 1]{Theorem}
\newtheorem{Theorem 8}[Theorem 1]{Theorem}
\newtheorem{Theorem 9}[Theorem 1]{Theorem}
\newtheorem{Theorem 10}[Theorem 1]{Theorem}

\title{\LARGE \bf
Privacy-Preserving Average Consensus via State Decomposition}

\author{Yongqiang Wang$^{1}$% <-this % stops a space
\thanks{ The work was supported in part by the National Science Foundation under Grants 1824014 and 1738902.}
\thanks{$^{1}$Yongqiang Wang is with the Department of Electrical and Computer Engineering, Clemson University, Clemson, SC 29634, USA
{\tt\small{yongqiw}@clemson.edu}
}%
}

\begin{document}

\maketitle
\thispagestyle{empty}
\pagestyle{empty}

%%%%%%%%%%%%%%%%%%%%%%%%%%%%%%%%%%%%%%%%%%%%%%%%%%%%%%%%%%%%%%%%%%%%%%%%%%%%%%%%
\begin{abstract}
Average consensus underpins key functionalities of distributed
systems ranging from distributed information fusion,
decision-making, distributed optimization, to load balancing and
decentralized control. Existing distributed average  consensus
algorithms require each node to exchange and disclose  state
information to its neighbors,  which is undesirable in cases where
the state is private or contains sensitive information. In this
paper, we propose a novel approach that avoids disclosing individual
state information in  average consensus  by letting each node
decompose its state into two sub-states. For each node, one of the
two sub-states   involves in computation and inter-node interactions
as if it were the original state, while the other sub-state
interacts only   with the first sub-state of the same node, being
completely invisible to other nodes. The initial values of the two
sub-states are chosen randomly but with their mean   fixed to the
initial value of the original state, which is key to guarantee
convergence to the desired consensus value. In direct contrast to
differential-privacy based privacy-preserving average-consensus
approaches which enable privacy by compromising accuracy in the
consensus value, the proposed approach can guarantee convergence to
the \emph{exact} desired value without any error. Not only is the
proposed approach able to prevent  the disclosure of a node's
initial state   to   honest-but-curious neighbors,  it can also
provide protection against inference by  external eavesdroppers able
to wiretap communication links.   Numerical simulations demonstrate
the effectiveness of the approach and   its advantages   over
state-of-the-art counterparts.
\end{abstract}

\section{Introduction}

As a building block of distributed computing, average consensus  has
been an active research topic in computer science and optimization
for decades \cite{Morris74,
    Lynch96}. In recent years, with the advances of wireless
communications and embedded systems, particularly the advent of
wireless sensor networks and the Internet-of-Things, average
consensus is finding increased applications in fields as diverse as
automatic control, signal processing, social sciences, robotics, and
optimization \cite{Olfati-Saber2007}.

Conventional average consensus approaches rely on the  exchange of
explicit state values among neighboring nodes to reach agreement on
distributed computation. Such a disclosure of state information has
two potential problems. First, it   breaches the privacy of
participating nodes  who may not want to disclose their state values
containing sensitive and private information. For example, a group
of individuals using average consensus to compute a common opinion
may want to keep each individual's opinion secret
\cite{citeulike84}. Another example is power systems where multiple
generators want to reach agreement on cost while keeping their
individual generation information private since the generation
information is sensitive in bidding the  right for  energy selling
\cite{fang2012smart}.  Secondly,   exchanging information through
wireless or wired communications is vulnerable to eavesdroppers
which try to steal information by tapping communication links. With
the increased number of reported attack events, preserving data
privacy has become an urgent need in many social and engineering
applications.
%, particularly in many real-time sensing
%and control systems such as power systems and wireless sensor
%networks.

To address the pressing need for privacy-preserving  average
consensus, one may resort to conventional secure multi-party
computation approaches such as Yao's Garbled Circuit
\cite{secure_computation1}, Shamir's Secret Sharing algorithm
\cite{secure_computation2}, or  many other recent advances
\cite{secure_computation3}. However, such general-purpose privacy
protecting approaches are both computationally and communicationally
too heavy for systems with fast-evolving behaviors particularly
cyber-physical systems which are subject to hard real-time
constraints. For example, Yao's Garbled Circuit has a computational
latency on the order of seconds \cite{kreuter2012billion}, whereas
the tolerable computational latency is on the order of milliseconds
for the real-time control of  connected automated vehicles
\cite{10ms2} and unmanned aerial vehicles \cite{latency_bound}.
Recently, several dedicated   privacy-preserving  solutions have
been proposed for average consensus \cite{Mo17,Erfan15,    Huang15,
Manitara13,nozari2015,katewa2017privacy,huang2012differentially}.
Most of these approaches rely on the idea of obfuscation to mask
true state values by injecting carefully-designed noise on the
states. One commonly used tool is differential privacy from the
database literature in computer science \cite{Erfan15,
Huang15,nozari2015,katewa2017privacy,huang2012differentially}.
However, obfuscation under differential privacy  affects  the
accuracy of average  consensus, preventing convergence to the exact
desired value. Another tool emerged recently is the correlated-noise
based obfuscation \cite{Manitara13,Mo17}, which can guarantee the
accuracy of average consensus.
 Observability based approaches have also been discussed in the dynamics and control community to
protect the privacy of multi-agent networks. The basic idea is to
design the interaction topology   to minimize the observability from
a certain node, which amounts to minimizing the node's ability to
infer the initial states of other  nodes in the network
\cite{kia2015dynamic,observability1,observability2}. However, both
the correlated-noise based   and the observability based approaches
are vulnerable to  adversary nodes which are directly connected to
   a target node as well as
all its neighbors \cite{ruan2019secure}.

 To improve
resilience to privacy attacks, another common approach is to employ
cryptography. However,  although cryptography based approaches can
easily enable privacy preservation with the assistance of  an
aggregator or third-party \cite{Gupta16}, like in cloud-based
control or computation
\cite{homomorphic_privacy1,homomorphic_privacy2,Lagendijk13}, their
extension to the completely \textit{decentralized} average consensus
problem \textit{in the absence of an aggregator or third-party} is
extremely hard due to the difficulties in  decentralized key
management. In fact, to our knowledge, except our recent result
\cite{ruan2017secure,ruan2019secure}, existing efforts
(\cite{Lazzeretti14,Freris16}) incorporate cryptography into
decentralized average consensus without giving participating nodes
access to the final consensus value (note that in \cite{Freris16}
individual participating nodes do not have access to the decryption
key to decrypt the final consensus value which is obtained in the
encrypted form, otherwise they will be able to decrypt intermediate
computations to access other nodes' states). Furthermore,
cryptography based approaches will also  significantly increase
communication and computation overhead (please see e.g.,
\cite{zhang2019admm} for detailed discussions), which is not
appropriate for systems with limited resources or systems with fast
evolving behaviors or subject to hard real-time constraints.

In this paper, we propose a state-decomposition based  approach that
can guarantee the privacy of all participating nodes in average
consensus without compromising accuracy. Our basic idea is to let
each node decompose its state into two sub-states with random
initial values. One sub-state succeeds the  role of the original
state in inter-node interactions while the other sub-state only
interacts with the first sub-state in the same node and thus is
completely invisible to outside nodes. To ensure  consensus to the
right average value, the initial values of the two sub-states are
randomly chosen but with  their mean fixed to the  initial   value
of the original state.   Different from existing
differential-privacy based
  approaches which  sacrifice accuracy for privacy, our approach can guarantee
convergence to the \textit{exact} average consensus value. Unlike
correlated-noise based  or observability  based approaches which
require a node to have at least one neighbor that is not directly
connected to the adversary to maintain privacy,  our approach can
guarantee privacy of a node even when the node and all its neighbors
are directly connected to the adversary. Furthermore, the approach
is completely decentralized and light-weight in computation, which
makes it easily applicable to resource-restricted systems. Numerical
simulation results are given to verify the results.

\section{Background}\label{sec:background}
\subsection{Average Consensus}
We first review the average consensus problem. Following the
convention  in \cite{Olfati-Saber2007}, we represent a network of
$M$ nodes as  a graph ${G=(V,\,E,\,\mathbf{A})}$ with node set
${V}=\{v_1,\,v_2,\,\cdots, v_M\}$, edge set ${E}\subset {V}\times
{V}$, and the adjacency matrix $\mathbf{A}=\big[a_{ij}[k]\big]$
denoting coupling weights which satisfy $a_{ij}[k]>0$ if
$(v_i,v_j)\in E$ and 0 otherwise. Here $k$ is time index, denoting
that $a_{ij}[k]$ could be time-varying. The set of neighbors of a
node $v_i$ is denoted as ${N}_i = \left\{v_j \in {V}| (v_i,v_j)\in
{E}\right\}$ and its cardinality is denoted as $|N_i|$.

Throughout this paper we make the following assumption:

\begin{Assumption 1}
  We assume that the graph is undirected and connected, i.e.,
$a_{ij}[k] = a_{ji}[k]$ holds for all $k\geq 0$ and there exists a
(multi-hop) path between any pair of nodes.
\end{Assumption 1}

 We represent the state variable of a node $i$ as $x_i[k]$. For the sake of simplicity, we assume scalar states. But as commented later in Remark 5,
 the results are easily extendable to the case where the state is a
 vector.
  To
achieve average consensus, namely convergence of all states $x_i[k]$
$(i=1,2,\cdots,M)$ to the average of initial values, i.e.,
$\frac{\sum_{i=1}^M x_i[0]}{M}$, the update rule    is   formulated
as \cite{olfati2007consensus}
\begin{equation}
\label{eq:dt} x_i[k+1] = x_i[k] + \varepsilon\sum_{v_j\in N_i}
a_{ij}[k] (x_j[k]- x_i[k])
\end{equation}
where $\varepsilon$ resides in the range $(0, \frac{1}{\Delta}]$
with $\Delta$ defined as
\begin{equation}\label{eq:Delta}
\Delta\triangleq\max_{i=1,2,\cdots,M}|N_i|
\end{equation}

It has been well known that  average consensus can be achieved if
the network is connected and there exists some $\eta>0$ such that
$\eta<a_{ij}[k]<1$ holds for all $k\geq 0$
\cite{nedic2010constrained}.

\subsection{Attack Model}
In the paper, we consider two types of adversaries:

\emph{ An honest-but-curious
    adversary}  is a node who
 follows all protocol steps correctly but is curious and collects
 received  data in an attempt to learn some information about other
 participating nodes.

\emph{An eavesdropper} is an external attacker   who knows the
network topology, and is able to wiretap communication links and
access exchanged messages.

Generally speaking, an eavesdropper is more disruptive than an
honest-but-curious node in terms of information breaches because it
can snoop messages exchanged on many channels whereas the latter can
only access the messages destined to it. However, an
honest-but-curious node does have one piece of information that is
unknown to an external eavesdropper, i.e., the internal initial
state  $x_i[0]$ is available if node $i$ is an honest-but-curious
node. We will systematically analyze the enabled privacy strength of
our approach against both adversaries.

\section{Privacy-Preserving Approach}\label{sec:proto}

\begin{figure*}[t!]
    \centering
    \includegraphics[width=.5\textwidth]{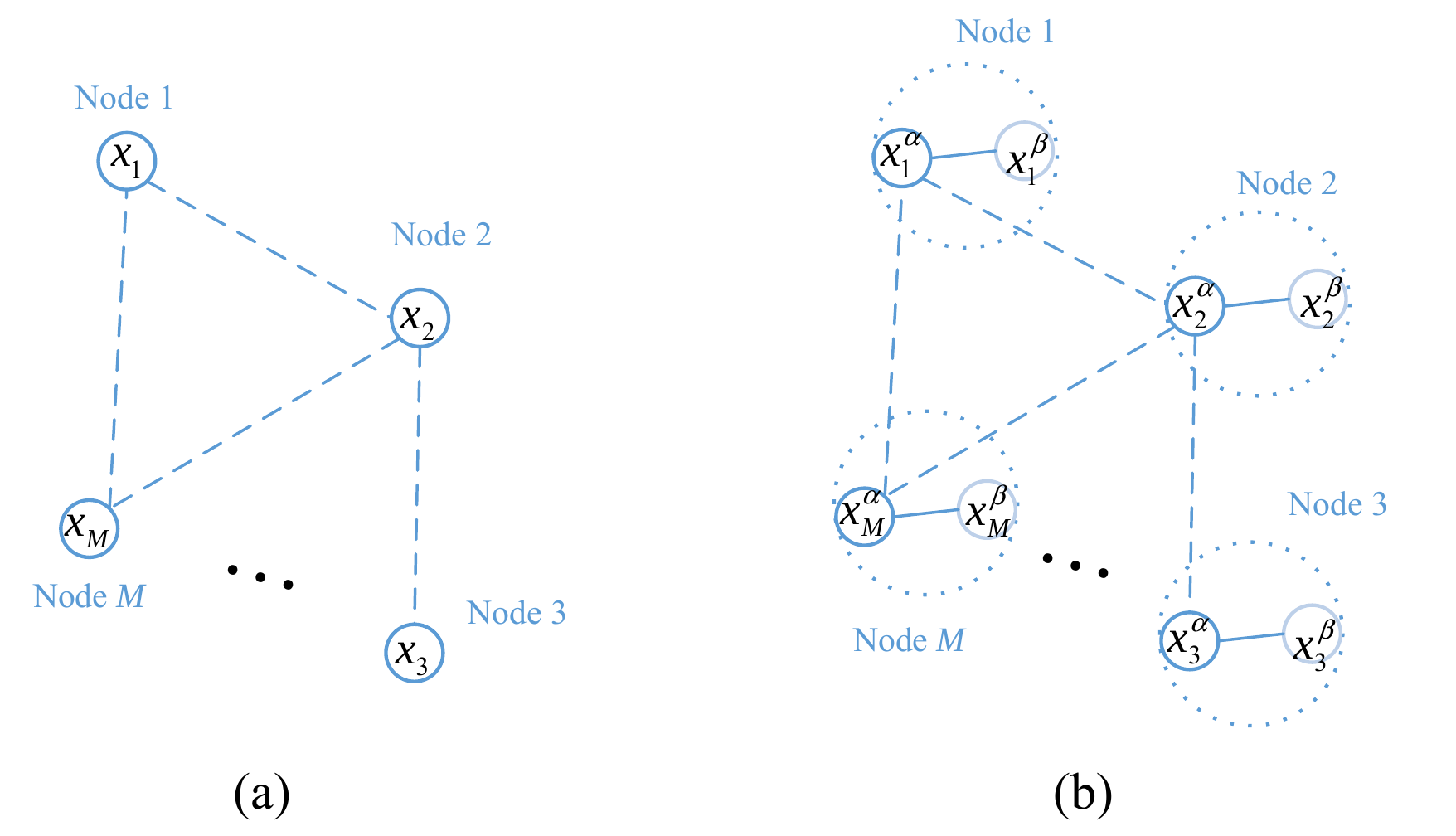}
        \vspace{-0.2cm}
    \caption{State-decomposition based privacy-preserving  average consensus. (a) Before state decomposition (b) After state decomposition}
    \label{fig:1}
\end{figure*}

The key idea of our approach is a decomposition mechanism:

\emph{Decomposition Mechanism:} We let each node decompose its state
$x_i$ into two sub-states $x_i^{\alpha}$ and $x_i^{\beta}$, with the
initial values $x_i^{\alpha}[0]$ and $x_i^{\beta}[0]$ randomly
chosen from the set of all real numbers under the constraint
$x_i^{\alpha}[0]+x_i^{\beta}[0]=2x_i[0]$ (cf. Fig. \ref{fig:1}). The
sub-state $x_i^{\alpha}$ succeeds the role of the original state
$x_i$ in  inter-node interactions and it is in fact the only state
value from node $i$ that can be seen by its neighbors. The other
sub-state $x_i^{\beta}$ also involves in the distributed interaction
by (and only by) interacting with $x_i^{\alpha}$. So the existence
of $x_i^{\beta}$ is invisible to neighboring nodes of node $i$,
although it directly affects the evolution of $x_i^{\alpha}$. Taking
node 1 in Fig. \ref{fig:1}(b) for example, $x_1^{\alpha}$ acts as if
it were $x_1$ in the inter-node interactions while $x_1^{\beta}$ is
invisible to nodes other than node $1$, although it affects the
evolution of $x_1^{\alpha}$. The coupling weight between the two
sub-states $x_i^{\alpha}$ and $x_i^{\beta}$ is symmetric and denoted
as $a_{i,\alpha\beta}[k]$. It is a design parameter and will be
elaborated later in the \emph{Weight Mechanism}.

Under the state-decomposition approach, the overall dynamics become
 \begin{equation} \label{eq:dt_decom}
 \left\{
 \begin{aligned}
 x_i^{\alpha}[k+1] &= x^{\alpha}_i[k] + \varepsilon\sum_{v_j\in N_i} a_{ij}[k](x^{\alpha}_j[k]- x^{\alpha}_i[k])\\
 &+\varepsilon a_{i,\alpha\beta}[k](x_i^{\beta}[k]-x_i^{\alpha}[k])\\
  x_i^{\beta}[k+1] &= x^{\beta}_i[k] + \varepsilon a_{i,\alpha\beta}[k](x_i^{\alpha}[k]-x_i^{\beta}[k])
 \end{aligned}
  \right.
 \end{equation}
subject to  $x_i^{\alpha}[0]+x_i^{\beta}[0]=2x_i[0]$.
\begin{Remark 1}
Compared with (\ref{eq:dt}), since every ``visible" sub-state's
number of neighbors is increased by 1, the upper bound on
$\varepsilon$ is reduced from $\frac{1}{\Delta}$ to
$\frac{1}{\Delta+1}$ with $\Delta$ defined in (\ref{eq:Delta}).
\end{Remark 1}

In conventional average consensus algorithms, the coupling weights
are required to be within $(0,\,1)$, which restricts the strength of
achievable privacy (as will be clear from the proof of Theorem 2).
We introduce the following weight mechanism to enable strong
privacy:

\emph{Weight Mechanism:} For $k=0$, we allow all weights $a_{ij}[0]$
and $a_{i,\alpha\beta}[0]$ to be arbitrarily chosen from the set of
all real numbers under the constraint $a_{ij}[0]=a_{ji}[0]$; For
$k=1,2,\cdots$, we require that there exists a scalar $0<\eta<1$
such that all nonzero $a_{ij}[k]$ satisfy $\eta\leq a_{ij}[k]<1$ and
all $a_{i,\alpha\beta}[k]$ satisfy $\eta\leq a_{i,\alpha\beta}[k]
<1$.

 In the
following, we first prove that under the approach, i.e., the
\emph{Decomposition Mechanism} and the \emph{Weight Mechanism}, all
states $x_i^{\alpha}$ and $x_i^{\beta}$ will converge to the same
average consensus value as in the conventional case (\ref{eq:dt}).
Then we rigorously analyze the privacy of participating nodes
enabled by the proposed approach in the presence of an eavesdropper
or honest-but-curious node.

\begin{Theorem 2}\label{th2}
        Under Assumption 1 and the \emph{Weight Mechanism},   all sub-states in (\ref{eq:dt_decom}) converge to the average consensus value of (\ref{eq:dt}), i.e.,
    \begin{equation}
    \lim_{k\to\infty} {x^{\alpha}_i}[k] =\lim_{k\to\infty}{x^{\beta}_i}[k]=    \frac{1}{M}\sum_{j=1}^M {x_j}[0]
    \end{equation}
\end{Theorem 2}

\noindent \textit{Proof}: Under the symmetric  weight assumption
$a_{ij}[k]=a_{ji}[k]$ and $a_{i,\alpha\beta}[k]$, one can easily
obtain that for the network after decomposition, the sum  of all
sub-states are always time-invariant. Therefore, even the weights
are allowed to be arbitrarily chosen from the set of all real
numbers for $k=0$, we always have
\begin{equation}\label{eq:sum}
 \frac{1}{2M}\sum_{j=1}^M ({x^{\alpha}_j}[0]+x^{\beta}_j[0])= \frac{1}{2M}\sum_{j=1}^M ({x^{\alpha}_j}[1]+x^{\beta}_j[1])
\end{equation}

Starting from $k=1$, as all coupling weights $a_{ij}[k]=a_{ji}[k]$
compose a connected graph, the \emph{Decomposition Mechanism} and
the \emph{Weight Mechanism} guarantee that all sub-states also
compose a connected graph. According to the result on average
consensus under time-varying weights \cite{nedic2010constrained},
average consensus can still be achieved, i.e., all sub-states
$x^{\alpha}_i$ and $x^{\beta}_i$ will converge to
$\frac{1}{2M}\sum_{j=1}^M ({x^{\alpha}_j}[1]+x^{\beta}_j[1])$, which
is equal to $\frac{1}{2M}\sum_{j=1}^M
({x^{\alpha}_j}[0]+x^{\beta}_j[0])$ according to (\ref{eq:sum}).
Further making use of the fact
$x_i^{\alpha}[0]+x_i^{\beta}[0]=2x_i[0]$  leads to the conclusion
that all sub-states converge to $\frac{1}{M}\sum_{j=1}^M {x_j}[0]$.
\hfill{$\blacksquare$}

\begin{Remark 3}
For the purpose of privacy-preservation, the  values of
$a_{i,\alpha\beta}[k]$ should be private to node $i$.
\end{Remark 3}

Next we rigorously analyze the enabled privacy against  an
honest-but-curious adversary or an external eavesdropping adversary.
To this end, we first give a definition of privacy.

\begin{Definition 1}
The privacy of the initial value $x_i[0]$ of node $i$ is preserved
if an adversary cannot estimate the value of $x_i[0]$ with any
guaranteed accuracy.
\end{Definition 1}

Definition 1 requires that an   adversary cannot even find a range
for a private value and thus is more stringent than the privacy
preservation definition considered in \cite{Mo17,Manitara13} which
defines privacy preservation as the inability of an adversary to
\emph{uniquely} determine the protected value. Next we show that
even by carefully observing a node's communication for multiple
steps, an adversary cannot infer the node's initial state with any
guaranteed accuracy.

\begin{Theorem 3}\label{theo:eavesdropper}
    Under the \emph{Decomposition Mechanism} and
the \emph{Weight Mechanism},  an
 honest-but-curious node $i$   cannot infer the initial state  $x_j[0]$ of node $j$ with  any guaranteed accuracy if node $j$ has at least one neighboring
 node $m$ who does not collude with node $i$ to infer  $x_j[0]$ (cf. Fig. \ref{fig:proof} for an illustrative example).
\end{Theorem 3}
\noindent \textit{Proof}:

\begin{figure}[t!]
    \centering
    \includegraphics[angle=0,width=.33\textwidth]{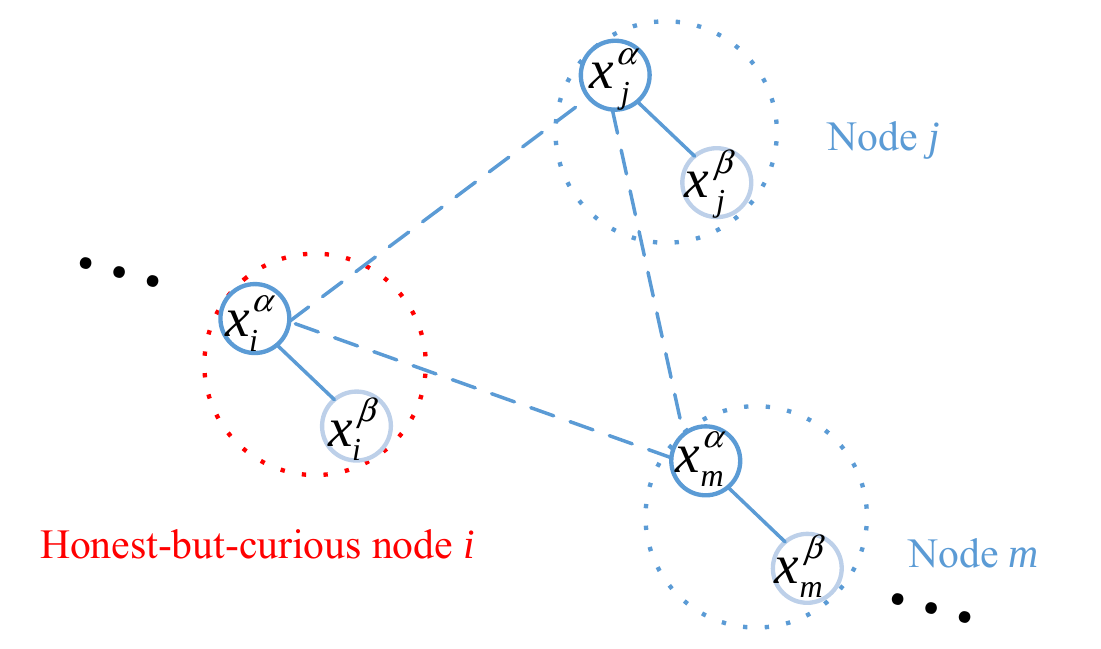}
    \caption{Connection configuration example for  Theorem 2.}
    \label{fig:proof}
\end{figure}

To prove that node $i$ cannot estimate $x_j[0]$ with any guaranteed
accuracy, we show that any arbitrary variation of $x_j[0]$ is
indistinguishable to node $i$, i.e., the information accessible to
node $i$ can be exactly the same even if $x_j[0]$ were changed to an
arbitrary value $\bar{x}_j[0]\neq x_j[0]$. We define the information
accessible to the honest-but-curious node $i$ at iteration $k$ as
$I_i[k]=\big\{a_{ip}[k]\big|_{v_p\in N_i},\,
x_{p}^{\alpha}[k]\big|_{v_p\in N_i}\,,
 x_i[k],\, x_i^{\alpha}[k],\,x_i^{\beta}[k],\,a_{i,\alpha\beta}[k]
 \big\}$. So as time evolves, the cumulated information accessible to node $i$ can be summarized as $I_i = \bigcup_{k=0}^{\infty} I_i[k]$.

To show that the privacy of the initial value $x_j[0]$ can be
preserved against node $i$, i.e., node $i$ cannot estimate the value
of $x_j[0]$ with any guaranteed accuracy, it suffices to show that
under any initial value $\bar{x}_j[0] \neq x_j[0]$ the information
accessible to node $i$, i.e., $\bar{I}_i$ could be exactly the same
as $I_i$, the cumulated information accessible to node $i$ under
$x_j[0]$. This is because the only information available for node
$i$ to infer the initial value $x_j[0]$ is $I_i$, and if $I_i$ could
be the outcome under any initial values of $x_j[0]$, then node $i$
has no way to even find a range for the initial value $x_j[0]$.
Therefore, we only need to prove that for any $\bar{x}_j[0] \neq
x_j[0]$, $\bar{I}_i = I_i$ could hold.

Next we show that there exist initial values of $x_m[0]$ and
coupling weights satisfying the requirements of the \emph{Weight
Mechanism} that make $\bar{I}_i = I_i$ hold under $\bar{x}_j[0] \neq
x_j[0]$. (Note that the alternative initial values of $x_m[0]$
should guarantee that the agents still converge to the original
average value after $x_j[0]$ is altered to $\bar{x}_j[0]$.) More
specifically, under the following initial condition
\begin{equation}\label{initial_values}
\begin{aligned}
\bar{x}_m[0] &= x_j[0] + x_m[0] -\bar{x}_j[0]\\
\bar{x}_j^{\alpha}[0]&=x_j^{\alpha}[0],\:\bar{x}_j^{\beta}[0]=2\bar{x}_j[0]-x_j^{\alpha}[0]\\
\bar{x}_m^{\alpha}[0]&=x_m^{\alpha}[0],\:\bar{x}_m^{\beta}[0]=2\bar{x}_m[0]-x_m^{\alpha}[0]\\
\bar{x}_q[0]&=x_q[0],\:\bar{x}_q^{\alpha}[0]=x_q^{\alpha}[0],\:\bar{x}_q^{\beta}[0]=x_q^{\beta}[0],\\
& \forall v_q \in \mathcal{V}\setminus \{v_j,v_m\}
\end{aligned}
\end{equation}
 and coupling weights
\begin{equation}\label{coupling_weight_bar}
\begin{aligned}
\bar{a}_{j,\alpha\beta}[0]&=\frac{x_j^{\beta}[0]-\bar{x}_j^{\beta}[0]+\varepsilon
{a}_{j,\alpha\beta}[0](x_{j}^{\alpha}[0]-x_{j}^{\beta}[0])
}{\varepsilon(\bar{x}^{\alpha}_j[0]-\bar{x}^{\beta}_j[0])}\\
\bar{a}_{m,\alpha\beta}[0]&=\frac{x_m^{\beta}[0]-\bar{x}_m^{\beta}[0]+\varepsilon
{a}_{m,\alpha\beta}[0](x_{m}^{\alpha}[0]-x_{m}^{\beta}[0])
}{\varepsilon(\bar{x}^{\alpha}_m[0]-\bar{x}^{\beta}_m[0])}\\
\bar{a}_{jm}[0]&=\frac{x_j^{\beta}[0]-\bar{x}_j^{\beta}[0]+\varepsilon
 a_{jm}[0](x_{m}^{\alpha}[0]-x_{j}^{\alpha}[0])
}{\varepsilon(\bar{x}^{\alpha}_m[0]- \bar{x}^{\alpha}_j[0])}\\
\bar{a}_{j,\alpha\beta}[k]&={a}_{j,\alpha\beta}[k],k=1,2\cdots\\
\bar{a}_{m,\alpha\beta}[k]&={a}_{m,\alpha\beta}[k],k=1,2\cdots\\
\bar{a}_{jm}[k]&= {a}_{jm}[k],k=1,2\cdots\\
\bar{a}_{q,\alpha\beta}[k]&={a}_{q,\alpha\beta}[k],\forall v_q \in \mathcal{V}\setminus \{v_j,v_m\},k=0,1,2\cdots\\
\bar{a}_{pq}[k]&=a_{pq}[k], \forall v_p,v_q \in
\mathcal{V},\{v_p,v_q\}\neq \{v_j,v_m\},\\
&\qquad   k=0,1,2\cdots
\end{aligned}
%a
\end{equation}
where ``$\setminus$" represents set subtraction, it can be easily
verified that $\bar{I}_i = I_i$ holds for any $\bar{x}_j[0] \neq
x_j[0]$. Note that the first equation in (\ref{initial_values}) is
used to guarantee that the consensus value does not change under the
alternative initial values $\bar{x}_j[0]$ and $\bar{x}_m[0]$.
Therefore, the honest-but-curious node $i$ cannot learn the initial
state of  node $j$ based on accessible information if node $j$ is
also connected to another node $m$ that does not  collude with node
$i$ to infer $x_j[0]$ (note that node $m$ is allowed to exchange
information with the honest-but-curious  node $i$ following the
protocol, as illustrated in Fig. \ref{fig:proof}).
\hfill{$\blacksquare$}

\begin{Remark 1}
In the derivation, the choice of coupling weights
$\bar{a}_{j,\alpha\beta}[0]$, $\bar{a}_{m,\alpha\beta}[0]$, and
$\bar{a}_{jm}[0]$ in (\ref{coupling_weight_bar}) guarantees that
starting from $k=1$, all sub-states under the alternative initial
value $\bar{x}_j[0]$ will be the same as those under the original
initial value $x_j[0]$, and hence all coupling weights can be the
same starting from $k=1$ under the two different initial value
conditions. Depending on the value of $\bar{x}_j[0]$, the weights
$\bar{a}_{j,\alpha\beta}[0]$, $\bar{a}_{m,\alpha\beta}[0]$, and
$\bar{a}_{jm}[0]$ could be outside the range $[\eta, \, 1)$, which
corroborates the necessity of allowing weights at $k=0$ to be
arbitrarily chosen from the  set of all real numbers in the
\emph{Weight Mechanism}. Note that since the sub-states
$\bar{x}_j^{\alpha}[0]= {x}_j^{\alpha}[0]$ and
$\bar{x}_m^{\alpha}[0]={x}_m^{\alpha}[0]$ are also arbitrarily
chosen from the  set of all real numbers, the possibility of them
making the denominators  in (7) equal to zero is negligible.
\end{Remark 1}

\begin{Remark 1}
    Our approach can protect the privacy of node $j$ even when node $j$ and all its neighbors are directly connected to the honest-but-curious neighbor $i$
 (cf. Fig. \ref{fig:proof}), which is not allowed in the privacy-preserving approaches in
\cite{Mo17,Manitara13}.
     This illustrates the advantage  of the proposed state-decomposition based
     approach.
\end{Remark 1}
%\begin{Remark 1}
%In the argument, since $c$ can be an arbitrary real number, the
%accessible information under a specific initial value of node 1 can
%be obtained after adding  an \emph{arbitrarily large offset} to the
%original initial value. So  the privacy-preserving approach can
%avoid an eavesdropper from estimating a finite range for the initial
%value. Therefore, the proposed approach can enable much stronger
%privacy protection than  existing studies in \cite{Mo17,Manitara13}
%which can only avoid an adversary from \emph{uniquely} determining
%the private value.
%\end{Remark 1}

Similar results can be obtained for the eavesdropping adversary
case:
\begin{Theorem 3}
   Under the \emph{Decomposition Mechanism} and
the \emph{Weight Mechanism}, an eavesdropper cannot infer the
initial state $x_j[0]$ of any  node $j$ with any guaranteed accuracy
   if node $j$ has at least one neighboring node $m$ whose interaction weight $a_{jm}[0]$ with node $j$  is inaccessible to the eavesdropper.
\end{Theorem 3}
\noindent \textit{Proof}: Following the line of reasoning in Theorem
2, we can obtain that any change in the initial value $x_j[0]$ can
be completely compensated by changes in $a_{jm}[0]$,
$a_{j,\alpha\beta}[0]$, and $a_{m,\alpha\beta}[0]$ that are
invisible to the eavesdropper. Therefore,  the accessible
information to the eavesdropper can be exactly the same even when
$x_j[0]$ were changed arbitrarily and hence the eavesdropper cannot
infer the initial value of node $j$ based on accessible information.
\hfill{$\blacksquare$}
\begin{Remark 2}
    The results are also applicable in the vector-state case. In fact, as long as the scalar state elements in the vector state have independent coupling weights, privacy can be naturally enabled in the vector-state case by applying results in this paper to individual scalar state
    elements.
\end{Remark 2}

\section{Numerical Comparison with Existing Results}\label{sec:example}

In this section, we numerically compare our  state-decomposition
based privacy-preserving approach with existing state-of-the-art
counterparts \cite{Mo17,Manitara13,huang2012differentially} to
confirm its   advantages.

For the convenience in comparison, we represent the internal state
of node $i$  as $x_i$ and its obfuscated version (used in state
exchange with neighbors) in existing obfuscation based approaches as
$x_i^+$. We considered a network of five nodes with interaction
topology and weights given in \cite{Mo17}, which, under our
formulation framework translate into $\varepsilon=\frac{1}{3}$ and
\begin{equation}
A = 0.75\begin{bmatrix}
0 & 1 & 0 & 0 & 1\\
1 & 0 & 1 & 0 & 0\\
0 & 1 & 0 & 0 & 1\\
0 & 0 & 0 & 0 & 1\\
1 & 0 & 1 & 1 & 0\\
\end{bmatrix}
\end{equation}

Without loss of generality, we suppose that an external eavesdropper
is interested in obtaining the initial   state of node 1 and
constructs the following observer to infer  $x_1[0]$:
    \begin{equation}
    \begin{aligned}
z[k+1]&=z[k]+x^+_1[k+1]\\
&\qquad-\left(x^+_1[k] + \varepsilon \sum_{v_j\in N_1}
a_{1j}(x^+_j[k]-x^+_1[k])\right)
    \end{aligned}
\label{eq:obs_update}
\end{equation}
where the initial observer state is set to $z[0] = x^+_1[0]$. We
assume that the eavesdropper has access to all weights except
$a_{12}[0]$ which was set to a random value 0.7  in the
 $k=0$ step of the observer.

\begin{figure}[tpb]
\hspace{-0.4cm}
    \includegraphics[width=0.55\textwidth]{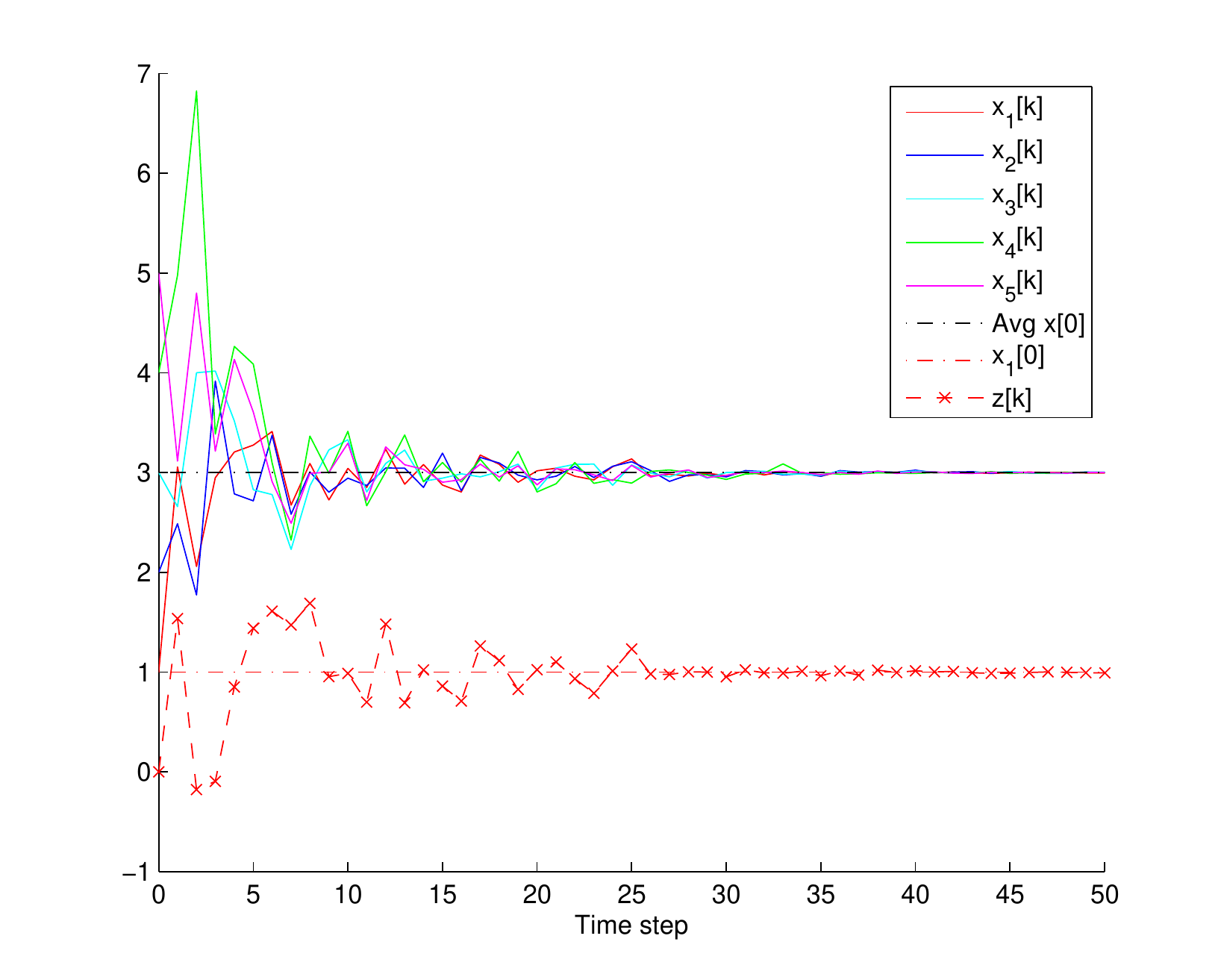}
    \caption{An eavesdropper can infer the initial state of node 1   under the privacy protocol in \cite{Mo17}.}
    \label{fig:mo_obs}
\end{figure}

Fig. \ref{fig:mo_obs} gives the evolution of the network states as
well as the eavesdropper's observer state   under the
privacy-preserving approach in \cite{Mo17}. The   initial states of
five nodes were set to $\left\{1,\,2,\,3,\,4,\,5\right\}$. It can be
seen that although convergence to the right average value is
achieved, the initial internal state of node 1, i.e., $x_1[0]$, can
also be inferred by the eavesdropper's observer $z[k]$.

\begin{figure}[tpb]
\hspace{-0.4cm}
    \includegraphics[width=0.55\textwidth]{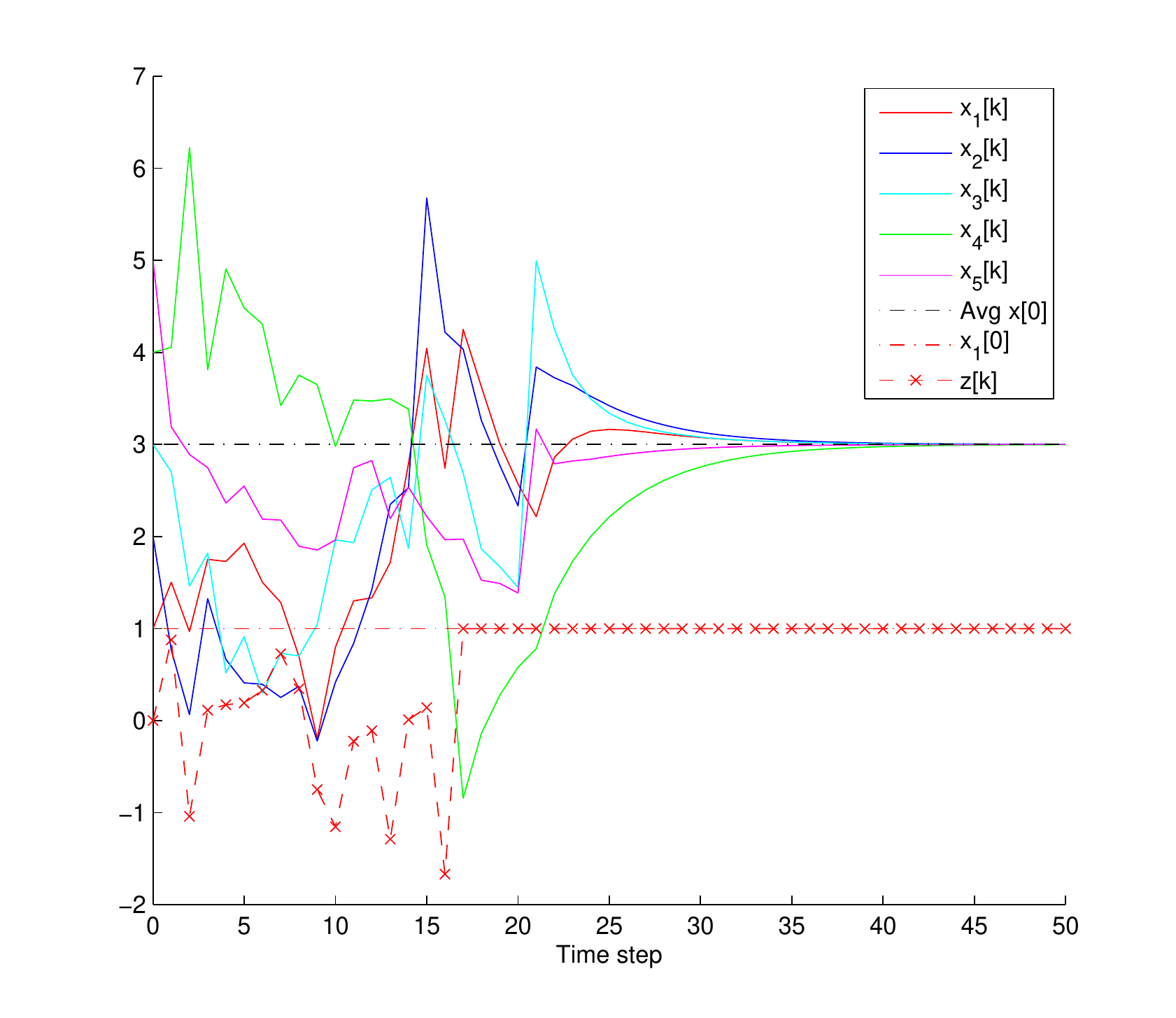}
    \caption{An eavesdropper can infer the initial state of node 1   under the privacy protocol in \cite{Manitara13}.}
    \label{fig:manitara_obs}
\end{figure}

Similar results were obtained using the approach in
\cite{Manitara13}, which   guaranteed accurate average consensus but
not the privacy of $x_1[0]$ (cf. Fig. \ref{fig:manitara_obs}).

\begin{figure}[tpb]
\hspace{-0.4cm}
    \includegraphics[width=0.55\textwidth]{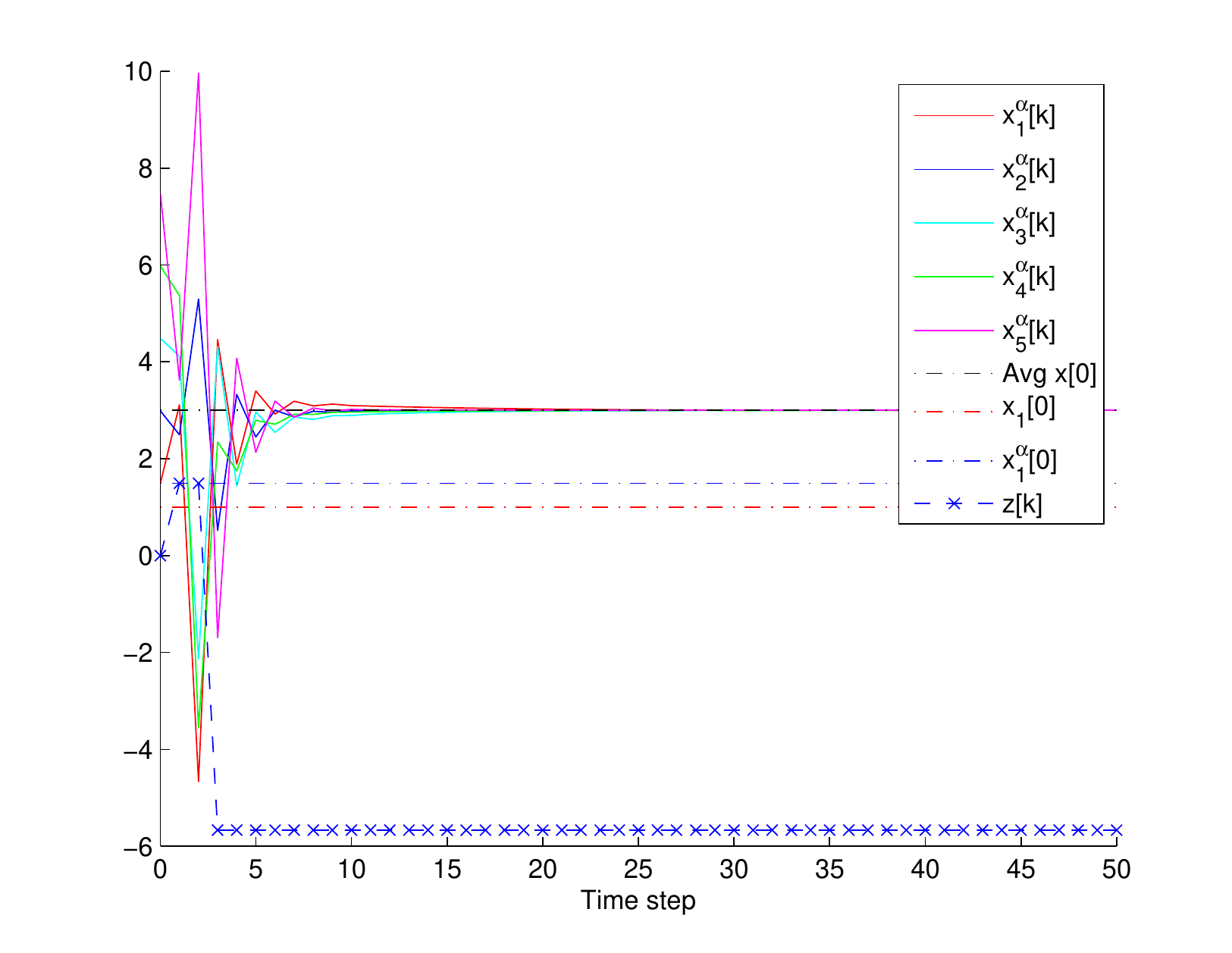}
    \caption{The proposed approach can achieve accurate average consensus while avoiding an eavesdropper from inferring the initial state $x_1[0]$ of node 1.}
    \label{fig:out_obs}
\end{figure}

 Based on the same setup, we also simulated the proposed state-decomposition based privacy-preserving
 approach. The coupling weights at $k=0$ were randomly chosen from
 $[-20,\,20]$.
 The results are given in Fig. \ref{fig:out_obs}, which confirms that the proposed  approach can protect the privacy of all nodes' initial values against
 an eavesdropper while achieving accurate average consensus.

\begin{figure}[tpb]
    \centering
    \includegraphics[width=0.45\textwidth]{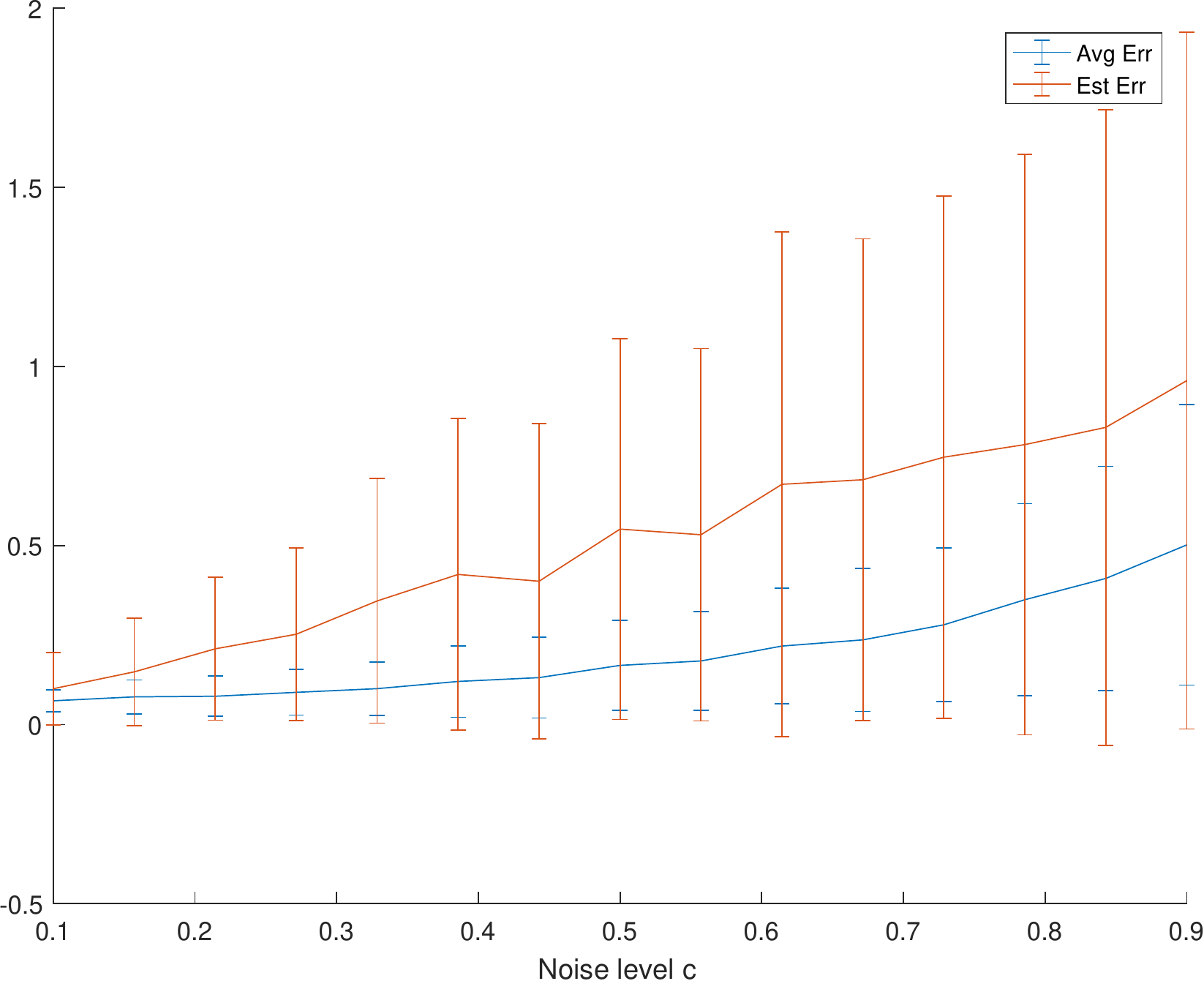}
    \caption{Although by using a high level of noise, the differential privacy approach in \cite{huang2012differentially} can avoid an eavesdropper from  accurately estimating the initial state of node 1
    (with large estimation error $\rm Est\: Err$), its   achievable  accuracy in average consensus also deteriorates  with an increased noise level, as reflected by a larger consensus error $\rm Avg\: Err$.}
    \label{fig:huang_Nfold}
\end{figure}

  It is worth noting that although differential-privacy based approaches such as \cite{huang2012differentially} can also protect the privacy of participating nodes' initial values, they also lead to errors in the final consensus value, as confirmed  by the simulation results in Fig. \ref{fig:huang_Nfold}. In such approaches, due to the trade-off between privacy and accuracy, when an application calls for higher  accuracy of the consensus result, the risk of disclosing one's initial state also becomes  higher.

%\begin{figure}[tpb]
%\hspace{-0.4cm}
%    \includegraphics[width=0.55\textwidth]{our_ob_enhanced.eps}
%    \caption{ Under random switching between sub-states, the proposed state-decomposition approach  can achieve accurate average consensus while avoiding an eavesdropper from inferring the initial state $x_1[0]$ of node 1.}
%    \label{fig:out_obs_enhanced}
%\end{figure}

\section{Conclusions}\label{sec:conclusion}
In this paper, we proposed a privacy-preserving approach for the
network average consensus problem based on state decomposition. In
contrast to differential-privacy based approaches which are subject
to a fundamental trade-off between enabled privacy and achievable
consensus accuracy, the approach is  able to enable privacy
preservation while guaranteeing accurate average consensus. It is
also superior to correlated-noise based obfuscation approaches which
can guarantee accurate average consensus but not resilience to
adversaries which are directly connected to a target node as well as
all its neighbors. Simulation results confirmed the theoretical
predictions.

\bibliographystyle{unsrt}

\bibliography{reference1}

\end{document}